\documentclass[11pt]{amsart}
\usepackage[top=3.5cm, bottom=3.5cm, left=2.5cm, right=2.5cm]{geometry}
\title{Random Walks on Vertex-Transitive Graphs with Moderate Growth}
\author{David Guo}

\usepackage{graphicx}
\usepackage{framed}
\usepackage{amsmath}
\usepackage{amssymb}
\usepackage{amsfonts}
\usepackage{mathrsfs}
\usepackage{array}
\usepackage{enumerate}
\usepackage{comment}
\usepackage{amsthm}  
\usepackage{mathtools}
\usepackage{physics}
\usepackage{kbordermatrix}
\usepackage[hidelinks]{hyperref}
\usepackage[capitalise]{cleveref}

\allowdisplaybreaks

\DeclareMathOperator{\diam}{diam}
\DeclareMathOperator{\Aut}{Aut}

\DeclareMathOperator{\Stab}{stab}
\DeclareMathOperator{\Cayley}{Cayley}

\providecommand{\R}{\mathbb{R}}
\providecommand{\N}{\mathbb{N}}

\renewcommand{\O}{\mathcal{O}}
\renewcommand{\a}{\alpha}
\renewcommand{\t}{\theta}

\providecommand{\pnnm}{ P_o^{2n+m}}

\providecommand{\pvnmx}{ P_v^{n+m}(x)}
\providecommand{\pvnmu}{ P_v^{n+m}(u)}
\providecommand{\ponv}{ P_o^{n}(v)}
\providecommand{\pon}{ P_o^{n}}

\providecommand{\agp}{ \Aut(\G,P)}

\providecommand{\ponm}{ P_o^{n+m}}
\providecommand{\ponmv}{ P_o^{n+m}(v)}

\providecommand{\e}{ \eta}

\providecommand{\g}{ \gamma}
\providecommand{\G}{ \Gamma}
\providecommand{\sv}{\sum_{v\in V}}

\providecommand{\su}{\sum_{u\in V}}

\providecommand{\h}{\frac{1}{2}}
\providecommand{\fc}{\frac{1}{C}}
\providecommand{\normi}[1]{\norm{#1}_\infty}

\theoremstyle{plain} 
\newtheorem{thm}{Theorem}[section] 
\newtheorem{lem}[thm]{Lemma} 
 
\newtheorem{cor}[thm]{Corollary}

\theoremstyle{definition} 

\newtheorem{eg}[thm]{Example}

\theoremstyle{remark} 
\newtheorem*{rem}{Remark}

\crefname{thm}{Theorem}{Theorems}
\Crefname{thm}{Theorem}{Theorems}

\crefname{lem}{Lemma}{Lemmas}
\Crefname{lem}{Lemma}{Lemmas}

\crefname{prop}{Proposition}{Propositions}
\Crefname{prop}{Proposition}{Propositions}

\crefname{cor}{Corollary}{Corollaries}
\Crefname{cor}{Corollary}{Corollaries}

\crefname{eg}{Example}{Examples}
\Crefname{eg}{Example}{Examples}

\DeclarePairedDelimiterX{\inp}[2]{\langle}{\rangle}{#1, #2}

\usepackage[numbers]{natbib}

\usepackage{url}

\begin{document}  
\date{}

\begin{abstract}
We study the mixing time of the uniform lazy random walk on a finite vertex-transitive graph, in which the probability of staying put equals the probability of moving to any given adjacent vertex. Diaconis and Saloff-Coste showed that, for Cayley graphs of moderate growth, the mixing time of this walk is quadratic in the diameter of the graph. Their main tools were bounds on the return probability and on the spectral gap of the walk. We extend their result to vertex-transitive graphs of moderate growth.
\end{abstract}

\maketitle  

\tableofcontents
\section{Introduction} \label{sec:intro}

Throughout this paper, we denote by $\G$ a finite, connected, undirected graph without multiple edges. We begin by introducing the definitions required for the statement of the main theorem. Recall that every random walk on $\G$ is associated with a stochastic matrix $P$ which contains the transition probabilities of the walk. For $u,v \in V$, we denote by $P_u$ the $u$-th row of $P$; this is a probability vector with non-negative entries that sum to one. The $(u,v)$-th entry of $P$ is denoted by $P(u,v)$ or $P_u(v)$, and represents the probability of moving from $u$ to $v$ given that we are at $u$. 

We define a \emph{transition system} $(\G, P)$ to be a graph $\G$ equipped with a transition matrix $P$ associated with a random walk on $\G$. We say that $(\G, P)$ is \emph{symmetric} if $P$ is symmetric. We denote by $U$ the matrix associated with the \emph{uniform distribution}, i.e.\ $U(u,v) = \frac{1}{N}$ for all $u, v \in V$. We write $\Aut(\G)$ for the group of graph automorphisms of $\G$, acting on the left. We define the set of automorphisms of $\G$ \ that \textit{preserve}
$P$ to be
\[
\Aut(\G,P) = \{ \alpha \in \Aut(\G) : \forall u,v \in V,\; P(u,v) = P(\alpha(u), \alpha(v)) \}.
\]
Observe that this is a subgroup of $\Aut(\G)$. We say that the system $(\G, P)$ is \emph{transitive} if $\Aut(\G,P)$ acts transitively on $\G$. In particular, if $(\G, P)$ is transitive, then $\G$ is vertex-transitive. Note that symmetry and transitivity of $(\G, P)$ do not imply one another.

Diaconis and Saloff-Coste \citep{MGA} showed that the uniform lazy random walk on a Cayley graph of moderate growth mixes in time quadratic in the diameter. Goel, Montenegro and Tetali \citep{MTB} later studied mixing times through the spectral profile, under a closely related notion of moderate growth. The aim of this paper is to show that the original methods of Diaconis and Saloff-Coste, namely bounds on the return probability and on the spectral gap, also work for general vertex-transitive graphs, with the same quadratic dependence on the diameter. Recall that, given constants $A$ and $d$, we say a finite connected vertex-transitive graph $\G$ has \emph{$(A,d)$-moderate growth} if
\[
\beta(j) \geq \frac{1}{A} \left( \frac{j}{\diam(\G)} \right)^d |\G| \quad \text{whenever} \quad 1 \leq j \leq \diam(\G).
\]

The main result of this paper is the following quadratic bound on the mixing time in $\ell^2$.

\begin{thm} \label{main_thm}
Let $(\G, P)$ be a symmetric transitive system. Suppose $\G$ is finite, connected, and has $(A,d)$-moderate growth. Let $\gamma = \diam(\G)$.

\begin{enumerate}[(i)]
\item \label{main_upper} For an upper bound, suppose that $\eta = \inf\{P(u,v) : u \sim v\} > 0$ and $\delta = \inf\{P(u,u) : u \in V\} \geq \frac{\eta}{2\gamma^2}$. Suppose further that $n$ is such that
\[
c := \frac{n\eta}{\gamma^2} - 4 > 0.
\]
Then
\[
\frac{\| (P^{n}-U)_o \|_{2}}{\| U_o \|_{2}} \leq B e^{-c},
\]
where $B = A^{1/2} 2^{1+\frac{d+d^2}{4}}$.

\item \label{main_lower} For a lower bound, suppose that $P(u,v)=0$ for any two non-adjacent vertices $u, v \in V$, and that $\gamma \geq A 2^{2d+2}$. Then
\[
\| (P^{n}-U)_o \|_{1} \geq e^{-c} \quad \text{for} \quad c = \frac{2^{4d+3} A^2 n}{\gamma^2}.
\]
\end{enumerate}
\end{thm}

Our proof follows the strategy of Diaconis and Saloff-Coste \citep{MGA}, who proved the analogue of \cref{main_thm} for Cayley graphs. We streamline their arguments, strengthen the results, and provide further clarification in some places.

The overall structure of our argument follows that of \citep{MGA}: \cref{P_plus}, \cref{bound_by_1_over_sqrt_2m}, \cref{bound_on_p2nm} and \cref{poly_growth_bound} correspond to Lemmas~2.4, 2.5, 2.6 and Theorem~2.3 of \citep{MGA} respectively; Section~5 generalises their Lemma~2.2; and the proof of \cref{main_thm} follows that of their Theorem~3.1. In each section below we indicate which part of their argument we are following. The main difference is that $\G$ need not be a Cayley graph, so we cannot use translations by group elements directly; instead we average over the stabiliser of a vertex in $\Aut(\G,P)$ (see \cref{change_sum_over_all_vetices_} and \cref{thm:bounds_on_Dirichlet_form}).

Recall that a symmetric transition matrix $P$ has real eigenvalues $1 = \pi_1 \geq \pi_2 \geq \cdots \geq \pi_N \geq -1$. Given such a matrix $P$, we set $\pi_* = \max\{|\pi_2|, |\pi_N|\}$, the second-largest eigenvalue of $P$ in absolute value. The following result motivates our study of the return probability and the eigenvalues of the random walk, forming our strategy for proving \cref{main_thm}.

\begin{lem} \label{bound_2_norm_by_e_vec}
Suppose $(\G, P)$ is a symmetric transitive system. Let $\pi_* = \max\{|\pi_2|, |\pi_N|\}$. Then for non-negative integers $n$ and $m$,
\begin{equation} \label{bound_aim}
\| (P^{n+m}-U)_o \|_2^2 \leq P^{2m}(o,o)\, \pi_*^{2n}.
\end{equation}
\end{lem}

\Cref{bound_2_norm_by_e_vec} is a direct generalisation of \citep[Lemma~2.1]{MGA}, which is stated in the setting of Cayley graphs. The structure of this paper is as follows. In Section 2, we introduce further notation and standard results required for this paper. In Section 3, we prove Lemma \ref{bound_2_norm_by_e_vec}. In Section 4, we bound the return probability $P^{2m}(o,o)$, and in Section 5, we bound the eigenvalue term $\pi_*^{2n}$. In Section 6, we combine these bounds to prove \cref{main_thm}; in fact, we also deduce mixing time bounds in any $\ell^p$ space from \cref{main_thm}. Finally, Tessera and Tointon  proved that vertex-transitive graphs with a suitable large diameter condition have moderate growth \citep{FST}, generalising work of Breuillard and Tointon and building on work of Breuillard, Green and Tao, Woess, Trofimov and Carolino (see the references in \citep{FST}); we show in Section 7 that our results imply that such graphs with large diameter also have quadratic mixing time, addressing a remark raised in the paper of Tessera and Tointon \citep[Remark 2.9]{FST}.

\subsection*{Acknowledgements}
I am very grateful to my PhD supervisor, Matthew Tointon, for his guidance and for reading drafts of this article.

\section{Preliminaries}

With the setup of \cref{sec:intro} in place, we introduce additional notation. Let $\G$ be a finite, connected, undirected graph without multiple edges. We write $V = V(\G)$ for the vertex set of $\G$. We often assume that $\G$ has $|V| = N$ vertices. For $u$, $v \in V$, we write $u \sim v$ if  $u$ and $v$ are connected by an edge; we denote by $d_\G(u,v)$ the number of edges in a shortest path from $u$ to $v$. Taking a starting vertex   $o \in V$,  we denote $|v| = d_\G(v,o)$ to be the number of edges in a shortest path from $o$ to $v$. Given $j \in \mathbb{N}$, we write $B_{\G}(v,j) = \{ u \in V : d(u,v) \leq j \} $ to be  the ball of radius $j$ centred at $v$.  Lastly, in the case when $\G$ is vertex-transitive, the size of the ball, $|B_{\G}(v,j)|$ is the same for every $v \in V$. We write $\beta_\G(j)$ for this cardinality. Also, we write $\gamma = \diam(\G)= \min\{j : \beta_\G(j)= |V|\}$, the diameter of the graph. We might often drop the subscript  $\G$ if it is clear from the context. Denote $\deg(\G)$ to be the degree of every vertex in $\G$.

The matrix $P$ associated with the \emph{simple random walk} is given by $P(u,v) = 1/\deg(\G)$ if $u \sim v$ and $P(u,v)=0$ otherwise; in particular, the walk never stays still. The \emph{lazy random walk} stays still with probability $1/2$ and otherwise moves to a uniformly random neighbour. We consider a random walk that differs from both of the above, which we call the \emph{uniform lazy random walk}: the probability of staying put is equal to the probability of moving to any given adjacent vertex, i.e.\ the walk moves from $u$ to a uniformly random vertex of $\{u\} \cup \{v : v \sim u\}$. This definition follows the approach of \citep{MGA}, who consider random walks on Cayley graphs whose generating sets contain the identity element. The transition matrix is
  \begin{equation} \label{simple_random_walk}
   P(u,v)=
    \begin{cases}
      \frac{1}{\deg(\G)+1} & \text{if }  u \sim v \text{ or } u=v , \\
      0 & \text{otherwise.}
    \end{cases}
  \end{equation}

As we mentioned in \cref{sec:intro}, a system $(\G, P)$ is transitive if the probability of going from any vertex to another is invariant under automorphisms.

For every $n \geq 1$ we have
\begin{equation} \label{inclusion_auto_group}
\Aut(\G,P) \subseteq  \Aut(\G,P^n). 
\end{equation}
In particular, if $(\G,P)$ is a transitive system, then so is $(\G,P^n)$ for every $n \geq 1$.

\begin{rem}
We will give an example when the other inclusion in (\ref{inclusion_auto_group}) does not hold. Take $\G = K_3$, the complete graph on $3$ vertices $v_1$, $v_2$ and $v_3$. Consider \renewcommand{\kbldelim}{(}
\renewcommand{\kbrdelim}{)}
\[
  P =\kbordermatrix{
    & v_1 & v_2 & v_3  \\
    v_1 & 1 & 0 & 0  \\
    v_2 & 0 & 0 & 1 \\
    v_3 & 0 & 1 & 0  \\
  }. 
\]
 Then $\Aut(\G,P^2) = \Aut(\G,I_3) \cong S_3$. Let $\a$ be an automorphism on $\G$ such that $\a(v_1) = v_2$. Then $P(v_1,v_1) =1\not = 0  =P(v_2,v_2)= P(\a(v_1),\a(v_1))$, so $\a \not \in \Aut(\G,P)$. 
\end{rem}

 Next, we introduce some standard definitions and results from functional analysis. Any real-valued function $f$ on the vertex set $V$ of $\G$, with $|V| = N$, can be viewed as a vector in $\mathbb{R}^N$.

Given a matrix $Q \in \mathbb{R}^{N \times N}$, we define the operator norm
\[
\| Q \|_{s \to t} = \sup_{\|f\|_s = 1} \|Qf\|_t.
\]
Recall that for $1 \leq s \leq t \leq \infty$, we also have
\begin{equation} \label{fun_norm_bound}
\| f \|_s \leq N^{\frac{1}{s} - \frac{1}{t}} \| f \|_t.
\end{equation}

We write $\langle f, g \rangle = \sum_{v \in V} f(v) g(v)$ for the inner product of real functions $f$, $g$ on $V$.
Let $P$ be a symmetric and stochastic matrix, the \emph{Dirichlet form} associated with $P$ on $V$ is defined by 
\[  \mathscr{E} _P = \mathscr{E} _P (f,f) = \inp{(I-P)f}{f}. \]
The following result gives us an expression of the Dirichlet form as a double sum over the vertices of $V$.
\begin{lem}  \label{change_u_v_sum_to_inner_product}
Let $P$ be a symmetric and stochastic matrix, for any real function $f$ on finite set $V$, we have 
\[  \mathscr{E} _P (f,f) = \frac{1}{2} \sum_{u,v \in V} (f(u)-f(v))^2 P(u,v). \]
\end{lem}

This identity is stated without proof in \citep[(2.7), p.~465]{NIF}. We include a short proof, which shows where each hypothesis is used.
\begin{proof}
Since $P$ is stochastic, $\sum_{v} P(u,v)=1$, so $\inp{(I-P)f}{f} = \sum_{u,v} f(u)^2 P(u,v) - \sum_{u,v} f(u)f(v)P(u,v)$. Since $P$ is symmetric, relabelling $u \leftrightarrow v$ gives $\sum_{u,v} f(u)^2 P(u,v) = \sum_{u,v} f(v)^2 P(u,v)$, so the first sum equals $\frac12\sum_{u,v}\big(f(u)^2+f(v)^2\big)P(u,v)$. Therefore $\inp{(I-P)f}{f} = \frac12\sum_{u,v}\big(f(u)^2+f(v)^2-2f(u)f(v)\big)P(u,v)$, which is the claim.
\end{proof}

\section{Bounding $\| (P^{n}-U)_o\|^2_2$ by the return probability and the eigenvalue}




We begin this section by establishing equality between the norm of a symmetric operator in the $l^2$ space and its eigenvalues. This is standard; see \cite[\S3.1, pp.~146--148]{HJ2}.

\begin{lem} \label{eigenvalue_of_symmetric}
Let $Q$ be a real symmetric $N \times N$ matrix with eigenvalues $\{\pi_i\}_{i=1}^N$, ordered so that $|\pi_1| = \max_i |\pi_i|$. Then
\[
\|Q\|_{2 \to 2} = |\pi_1|.
\]
\end{lem}

Recall that a symmetric transition matrix $P$ has real eigenvalues $1 = \pi_1 \geq \pi_2 \geq \cdots \geq \pi_N \geq -1$.  The above lemma applies to our case as follows. 

\begin{lem} \label{eignvalue_of_p-U}
Let $P$ be a real symmetric transition matrix. Let $1 = \pi_1 \geq \pi_2 \geq \cdots \geq \pi_N \geq -1$ be its eigenvalues, with corresponding orthogonal eigenvectors $f_1 = (1,\ldots,1), f_2, \ldots, f_N$. Then:
\begin{enumerate}[(i)]
\item \label{eignvalue_of_U}
$f_1, f_2, \ldots, f_N$ are also eigenvectors of $U$, corresponding to the eigenvalues $1, 0, \ldots, 0$, respectively.
\item \label{operator_2_2_bound}
$\|P - U\|_{2 \to 2} = \pi_*$, where $\pi_* = \max\{|\pi_2|, |\pi_N|\}$.
\end{enumerate}
\end{lem}
Part (i) is immediate, and part (ii) is stated in the proof of \citep[Lemma~2.1]{MGA}.

For a transitive system $(\G, P)$, the entries of any two rows are permutations of each other: if $u,v \in V$ and $\alpha \in \Aut(\G,P)$ satisfies $\alpha(u)=v$, then $P(u,w) = P(v,\alpha(w))$ for all $w \in V$. The following result relates the operator norm to the row vectors of such matrices.
\begin{lem}\label{holder_conjugate_norm_equality}
Let $T$ be an $N \times N$ matrix with rows and columns indexed by a finite set $V$. Suppose that $T$ satisfies
\begin{equation} \label{permutation_property}
\forall u,v \in V,\; \exists \sigma \in \Aut(V) \text{ such that } \forall x \in V,\; T_u(x) = T_v(\sigma(x)).
\end{equation}
Let $1 \leq s,t \leq \infty$ with $\frac{1}{s} + \frac{1}{t} = 1$. Then, for any $o \in V$,
\begin{enumerate}[(i)]
\item \label{ss}
if $T$ is additionally symmetric, $\|T\|_{s \to s} \leq \|T_o\|_1$;
\item \label{s_infty}
$\|T\|_{s \to \infty} = \|T_o\|_t$.
\end{enumerate}
\end{lem}
It is clear that any power of the transition matrix of a symmetric transitive system satisfies \eqref{permutation_property} and is itself symmetric.

\begin{proof}
By \eqref{permutation_property}, for every $u\in V$ there is a bijection $\sigma$ of $V$ with $T_u=T_o\circ\sigma$, so $\|T_u\|_p=\|T_o\|_p$ for all $u$ and all $p\in[1,\infty]$.

(i) By \eqref{permutation_property} and symmetry, every row and every column of $T$ has $\ell^1$ norm $\|T_o\|_1$. The Schur test then gives $\|T\|_{s\to s}\le\|T_o\|_1$ for all $1\le s\le\infty$.

(ii) For $\|f\|_s=1$, H\"older's inequality gives $|(Tf)(u)|\le \|T_u\|_t\|f\|_s=\|T_o\|_t$, so $\|T\|_{s\to\infty}\le\|T_o\|_t$. Equality holds at $u=o$ for $f=h/\|h\|_s$, where $h=|T_o|^{t-1}\operatorname{sgn}(T_o)$ if $1<t<\infty$, $h=\operatorname{sgn}(T_o)$ if $t=1$, and $h=\delta_y$ with $y$ a maximiser of $|T_o|$ if $t=\infty$.
\end{proof}

We take note of these two special cases of the above lemma: 
\begin{enumerate} [(i)]
\item  
If $P$ is a symmetric stochastic matrix that satisfies property (\ref{permutation_property}), then for $ 1 \leq s \leq \infty$,  
\begin{equation} \label{ssnorm_leq_1_norm}
\|P\|_{s\rightarrow s} \leq 1 
\end{equation}
\item
When we take $s = t = 2$ in Lemma \ref{holder_conjugate_norm_equality} (\ref{s_infty}), we have
\begin{equation} \label{2_infinity_operator_norm_2_norm}
\|T\|_{2\rightarrow \infty} = \|T_o\|_2. 
\end{equation}

\end{enumerate}

The following result is useful when we compute the difference between the $n$-step transition matrix and the uniform matrix. 

\begin{lem} \label{pu_up_u}
\begin{enumerate}[(i)]
\item  For any symmetric stochastic matrix  $ P$, we have $PU = UP = U$. 
\item Furthermore, for any $n \ge 1$, we have $P^n - U = (P-U)^n$. 
\end{enumerate}
\end{lem}

\begin{proof}
Since $P$ is symmetric and stochastic, it is doubly stochastic, so $PU = UP = U$. Also $U^2 = U$. Hence $(P-U)^n = P^n - U$ for $n \geq 1$.
\end{proof}

In the following lemma, we establish an equality between the return probability and the norm of the power of the transition matrix. Following from this, we also take a look at the decreasing pattern in the sequence $\{P^{2m}(o,o)\}_{m=1}^\infty$. 

\begin{lem} \label{some_bounds}
Suppose $(\G,P)$ is a symmetric transitive system, then for any $m \geq 0$, we have 
\begin{enumerate}[(i)] 
\item \label{2m_power_equal_2norm_squared_equal_inf_norm}    $P^{2m}(o,o) = \|P^m_o \|^2_2 = \|P^{2m}_o \|_\infty$;

\item \label{2m+2_power_bdd_by_2m_power}  $\|P^{m+1}_o \|_\infty \leq \|P^{m}_o \|_\infty$. Therefore, $P^{2m+2}(o,o) \leq P^{2m}(o,o)$. 
\item  $P^{2m+1}(o,o) \leq P^{2m}(o,o)$.

\end{enumerate}
\end{lem}
Part (i) is stated in the proof of \citep[Lemma~2.1]{MGA}. Part (ii) follows from H\"older's inequality: $\|P_o^{m+1}\|_\infty\le\|P_o^m\|_\infty\max_v\|P_v\|_1=\|P_o^m\|_\infty$. Part (iii) is stated in \citep[Lemma~2.4]{MGA}.

\begin{eg}
In general, we do not have  $P^{m+1}(o,o) \leq P^{m}(o,o)$ for $ m \geq 0$. For example, take $P= 
\big(\begin{smallmatrix}
  0 & 1\\
  1 & 0
\end{smallmatrix}\big)$, we have   $P^2= 
\big(\begin{smallmatrix}
  1 & 0\\
  0 & 1
\end{smallmatrix}\big)$.
\end{eg}

 \begin{rem} \label{rem:some_bounds_general}
Note Lemma \ref{some_bounds} \eqref{2m_power_equal_2norm_squared_equal_inf_norm} extends, in absolute-value form, to matrices that need not be stochastic or non-negative. Suppose $T$ is a real symmetric $N \times N$ matrix indexed by $V$ satisfying the permutation property \eqref{permutation_property} (with respect to some subgroup of $\Aut(\G)$ preserving $T$). As for \eqref{inclusion_auto_group}, by induction on $m$, every power $T^m$ also satisfies \eqref{permutation_property}, witnessed by the same automorphisms. Consequently, for every $m \geq 0$,
\begin{equation} \label{signed_2m_bound}
\|T^{2m}_o\|_\infty := \max_{v \in V} |T^{2m}(o,v)| = T^{2m}(o,o) = \|T^m_o\|_2^2.
\end{equation}
Indeed, symmetry gives
\[
T^{2m}(o,o) = \sum_{u \in V} T^m(o,u)^2 = \|T^m_o\|_2^2 \geq 0,
\]
and for any $v \in V$, Cauchy--Schwarz together with $\|T^m_v\|_2 = \|T^m_o\|_2$ (from the permutation property) gives
\[
|T^{2m}(o,v)| \leq \|T^m_o\|_2 \|T^m_v\|_2 = \|T^m_o\|_2^2 = T^{2m}(o,o),
\]
with equality at $v=o$; taking the maximum over $v$ gives \eqref{signed_2m_bound}.

We apply this to $T = P-U$. Since $U(u,v)=1/N$ is constant, every $\alpha \in \Aut(\G,P)$ automatically preserves $U$, so $\Aut(\G,P) \subseteq \Aut(\G,P-U)$; as $(\G,P)$ is transitive, so is $(\G,P-U)$. By Lemma \ref{pu_up_u}, $(P-U)^m = P^m-U$, so \eqref{signed_2m_bound} applies to $T=P-U$ for every $m \geq 0$.
\end{rem}


We can now give the proof for the result which we stated in the introduction. The proof below is the argument of \citep[Lemma~2.1]{MGA}, adapted to vertex-transitive graphs.


\begin{proof} [Proof of Lemma \ref{bound_2_norm_by_e_vec}]

For any real function $f$ on $V$, we have
\begin{align} \label{norm_bound}
 \|(P^{n+m} - U)f \| _\infty &= \| P^m(P^n-U)f \|_\infty  & \text{by Lemma \ref{pu_up_u}} \nonumber\\ 
							&\leq \| P^m \|_{2 \rightarrow \infty} \norm{P^n-U}_{2\rightarrow 2} \norm{f}_2. 						  
\end{align}

It follows that,
\begin{align*}
\| (P^{n+m}-U)_o\|^2_2 &= \norm{P^{n+m}-U}^2_{2\rightarrow \infty} &\text{by property (\ref{permutation_property}) for $T=P^{n+m}-U$ and equality (\ref{2_infinity_operator_norm_2_norm}}) \\
						& \leq \norm{P^m}_{2\rightarrow \infty} ^2 \norm{P^n-U}_{2\rightarrow 2} ^2 &\text{by  (\ref{norm_bound})  } \\
						& =  \| P^m_o\|^2_2  \norm{(P-U)^n}_{2\rightarrow 2}^2  &\text{by equality (\ref{2_infinity_operator_norm_2_norm}) and Lemma \ref{pu_up_u}  } \\
						& \leq P^{2m}(o,o) \pi_*^{2n} &\text{by Lemma \ref{some_bounds} (\ref{2m_power_equal_2norm_squared_equal_inf_norm}) and \ref{eignvalue_of_p-U} (\ref{operator_2_2_bound})}
\end{align*}
which completes the proof. 
\end{proof}

\section{The decay of the sequence $P^{2m}(o,o)$ }

We will be using the following standard results from linear algebra:

\begin{lem} \label{self_adjoint}
\begin{enumerate}[(i)]
\item \label{positive-semi-definite-iff-eigenvalue-are-non-negative}
Let $P$ be a symmetric matrix. Then $P$ is positive semi-definite if and only if all its eigenvalues are non-negative. 
\item \label{self_adjoint_root}
A positive semi-definite matrix $T$ has a unique positive semi-definite square root $T^{1/2}$. In particular, $T^{1/2}$ is symmetric. 
\item \label{symmetric_is_self_adjoint}
A linear map on $\R^N$ defined by a symmetric matrix is self-adjoint. 
\end{enumerate}
\end{lem}
When we have a semi-definite matrix $T$, we write $T^{1/2}$ to be its unique positive semi-definite square root. 
\begin{eg} \label{I-P_example}
Denote $I$ as the $N$ by $N$ identity matrix. Recall that the eigenvalues of a stochastic symmetric matrix $P$ are bounded above by $1$. Therefore, $I - P$ is positive semi-definite by Lemma \ref{self_adjoint} (\ref{positive-semi-definite-iff-eigenvalue-are-non-negative}). Hence, by Lemma \ref{self_adjoint} (\ref{self_adjoint_root}) and (\ref{symmetric_is_self_adjoint}), $(I - P)^{1/2}$ is self-adjoint. 
\end{eg}  
Let $(\G, P)$ be a transitive system. We will write $P^+ = (I+P)/2$, which is still a stochastic matrix. We define  $\eta=\eta_P=\inf\{P(u,v):u\sim v\}$ and similarly  $\eta^+=\eta^+_P=\inf\{P^+(u,v):u\sim v\}$. 

We now examine some properties of the variables associated with the matrix $P^+$. \Cref{P_plus} generalises \citep[Lemma~2.4]{MGA} from random walks on finite groups to symmetric transitive systems.

\begin{lem} \label{P_plus}
Let $(\G, P)$ be a symmetric transitive system. Then
\begin{enumerate}[(i)]
\item 
$P^+$ is positive semi-definite;

\item 
$\eta^+ = \frac{1}{2} \eta$;
\item
For $m \geq 0$, we have $P^{2m}(o,o) \leq   2  (P^+)^{2m}(o,o)$.

\end{enumerate}

\end{lem}

\begin{proof}
\begin{enumerate}[(i)]

\item 
Similar to \cref{I-P_example}, we know that the eigenvalues of $P$ are bounded below by $-1$, so the eigenvalues of $P^+$ are bounded below by 0, the result follows from Lemma \ref{self_adjoint} (\ref{positive-semi-definite-iff-eigenvalue-are-non-negative}). 
\item 
This follows from the fact that the off-diagonal elements of $I$ are $0$, i.e. $I(u,v) = 0$ for $u\sim v$. 
\item

We have 
\begin{align*}
(P^+)^{2m}(o,o) &= \left(\frac{I +P}{2}\right)^{2m}(o,o)\\
				&= \frac{1}{2^{2m}} \sum_{j=0}^{2m} \binom {2m}  {j} I^{2m-j} P^j(o,o) &\text{using Binomial expansion} \\
				& \geq \frac{1}{2^{2m}} \sum_{j=0}^{m} \binom{2m }{2j} P^{2j}(o,o)  \\
				&	\geq	 	P^{2m}(o,o) \left( \frac{1}{2^{2m}}  \sum_{j=0}^{m} \binom{2m }{ 2j} \right) & \text{ by Lemma \ref{some_bounds} (\ref{2m+2_power_bdd_by_2m_power}}) \\
				&= \frac{P^{2m}(o,o)}{2}.
\end{align*}
\end{enumerate}
\end{proof}

 The proof of the next lemma requires the following observation from calculus. Given $m \in \N$, consider the function $h:[0,1] \rightarrow \R$, $h(x)=(1-x)x^{2m}$, by calculus, $h$ obtains its maximum at $x = \frac{2m}{2m+1}$. 

Recall a standard analysis result that the sequence $ \left\{ \left( \frac{n}{n+1} \right)^n \right\}_{n=1}^\infty$ is a positive decreasing sequence. Therefore,\begin{equation} \label{bound_on_h}
\sup_{0 \leq x \leq 1} h(x)  = \left(\frac{1}{2m+1} \right) \left( \frac{2m}{2m+1} \right) ^{2m} \leq  \left( \frac{1}{2m+1} \right) \left(  \frac{4}{9} \right) \leq \frac{1}{4m}. 
\end{equation}

\begin{lem} \label{bound_by_1_over_sqrt_2m}
Suppose $P$ is an $N$ by $N$  positive semi-definite stochastic matrix (hence symmetric) and $m \geq 1$. Then \[\|(I-P)^\h P^m \|_{2 \rightarrow 2} \leq \frac{1}{2\sqrt{m}}. \]
\end{lem}

\begin{proof}
This is \citep[Lemma~2.5]{MGA} (due to Coulhon and Saloff-Coste). They state it for random walks on finite groups, but their proof applies verbatim to any positive semi-definite stochastic matrix: writing $f=\sum_i a_i f_i$ in an orthonormal eigenbasis of $P$, with eigenvalues $\pi_i\in[0,1]$, we have
\[
\|(I-P)^{1/2}P^m f\|_2^2=\sum_i (1-\pi_i)\pi_i^{2m}a_i^2\le \frac{1}{4m}\|f\|_2^2
\]
by \eqref{bound_on_h}, and the result follows by taking square roots.
\end{proof}

The following lemma is an orbit--stabiliser count: for a group $H$ acting transitively on $V$, summing $f(\alpha(v))$ over $\alpha\in H$ gives the size of a stabiliser times the sum of $f$ over $V$.

\begin{lem} \label{transitive_action}
Let $V$ be a finite set and $H$ a group acting on $V$. Let $f$ be a real function on $V$ and let $v$ be an element in $ V$.
Consider the following sums:

\begin{enumerate}[(i)]
\item 
$S_1=\sum_{\a \in H}f(\a(v))$
\item
$S_2= \su \sum_{\substack{\a \in H,\\  \alpha(v)=u}} f(\a(v))$
\item
$S_3= \su \sum_{\substack{\a \in H,\\  \alpha(u)=v}} f(\a(v))$.
\end{enumerate}

Then $S_1 = S_2=S_3$ for any group action $H$ on $V$. In the case where  $H$ acts transitively on $V$, we have $ |\Stab(w)|$, the size of the stabiliser of $w$ under  $H$,  is the same for every $w$ in $V$, we denote by $C$ to be this cardinality. Then the following sum also coincides with the above ones:
\begin{enumerate}[(iv)]
\item $S_4=C \su f(u)$.
\end{enumerate}

\end{lem}

\begin{proof}

($S_1 = S_2=S_3$) This follows from the following observation:
\[\cup_{u \in V } \{ \a \in H | \a(v)= u \} = H = \cup_{u \in V}  \{ \a \in H | \a(u)= v \}. \]
($S_2=S_4$ when the action is transitive) Let $ u \in V$, then there exists $\beta \in H$, such that $\beta(v)=u$. Also, $ \{ \a \in H | \a(v)= u \} = \beta\Stab(v) $. Hence $|\{ \a \in H | \a(v)= u \}| = |\Stab(v)| = C$. 
\end{proof}



We will see how the above lemma can be applied to a transitive system.  The next lemma is where our argument departs from that of \citep{MGA}: for a Cayley graph the relevant comparison can be made using translations by group elements (see the first displayed formula in the proof of \citep[Lemma~2.6, p.~8]{MGA}, where the value at $x$ is compared with the value at $xz$ for a generator $z$), whereas for a general vertex-transitive graph we average over the stabiliser of a vertex, using \cref{transitive_action}.

\begin{lem} \label{change_sum_over_all_vetices_}
Let $(\G, P)$ be a symmetric transitive system and fix $o\in V$. Let $T$ be an $N\times N$ matrix indexed by $V$ that is invariant under $\agp$, i.e. $T(\a(v),\a(x))=T(v,x)$ for all $\a\in\agp$ and $v,x\in V$. Then for all $x\sim u$,
\[ \sv  \big ( T_v(x) - T_v(u) \big )^2 \e  \leq \sv \sum_{w \in V} \big ( T_o(w) - T_o(v) \big )^2 P(w,v), \] where $\eta=\inf\{P(u,v):u\sim v\}$.
\end{lem}
By \eqref{inclusion_auto_group}, every power $P^k$ with $k\geq 1$ satisfies the invariance hypothesis, so \cref{change_sum_over_all_vetices_} applies to $T=P^{n+m}$; this is how it is used in the proof of \cref{bound_on_p2nm}.
\begin{proof}
Let $C = |\Stab(v)|$, the size of the stabiliser of $v$ in $\agp$. We have
\begin{align*}
 & \quad \sv  \bigg [ T_v(x) - T_v(u) \bigg ]^2 \e  \\
 &= \sv  \fc \sum_{\substack{\a \in \agp,\\  \alpha(v)=o}} \bigg [ T_v(x) - T_v(u) \bigg ]^2 \e   \\
 &= \fc \sv   \sum_{\substack{\a \in \agp,\\  \alpha(v)=o}} \bigg [ T_{\a(v)}(\a(x))-  T_{\a(v)}(\a(u)) \bigg ]^2 \e \quad \text{since $\a\in\agp$ preserves $T$} \\
 &= \fc \sv   \sum_{\substack{\a \in \agp,\\  \alpha(v)=o}} \bigg [ T_o(\a(x))-  T_o(\a(u)) \bigg ]^2 \e \\
  &= \fc \sum_{\a \in \agp}  \bigg [ T_o(\a(x))-  T_o(\a(u)) \bigg ]^2 \e  \qquad \text{by Lemma \ref{transitive_action}} (S_1=S_3)\\
 &\leq  \fc \sum_{\a \in \agp}  \sum_{\substack{{w \in V},\\  w \sim \a(x)} }   \bigg [ T_o(\a(x))-  T_o(w) \bigg ]^2 \e \qquad  \text{since $\a(x) \sim \a(u)$} \\
 &\leq   \fc   C \sv  \sum_{\substack{{w \in V},\\  w \sim v}}  \bigg [ T_o(v)-  T_o(w) \bigg ]^2  P(v,w) \qquad \text{by Lemma \ref{transitive_action}} (S_1=S_4) \\
 &\leq   \sv \sum_{w \in V} \big ( T_o(w) - T_o(v) \big )^2 P(w,v).\qedhere
 \end{align*} 
\end{proof}
The following lemma establishes a connection between the rate of decrease of the return probability $P^{2m}$   and the volume growth of $\G$. 
\begin{lem} \label{bound_on_p2nm}
Let $(\G, P)$ be a symmetric transitive system with $P$ being a positive semi-definite matrix and $\eta=\inf\{P(u,v):u\sim v\}>0$.  Then, $\forall n,m \geq 0$
\[
P^{2n+m}(o,o) \leq \frac{2}{\beta(r(n,m))}, \quad \text{where} \quad r(n,m)= \left(\frac{\eta}{2} \right)^{\frac{1}{2}} m^{\frac{1}{2}} \frac{P^{2n+m}(o,o)}{P^{2n}(o,o)}.
\]
\end{lem}

Apart from the use of \cref{change_sum_over_all_vetices_}, the proof below follows that of \citep[Lemma~2.6]{MGA}.

\begin{proof}

Let $x \in V$ and $u \sim x$.  Then
\begin{align*}
    & \quad |\pnnm(x)-\pnnm(u)|\\
  &= \sv \big | \ponv \pvnmx - \ponv  \pvnmu \big | \\
  &= \sv \big | \ponv \big | \big | \pvnmx - \pvnmu \big | \\
  &\leq  \left( \sv \big | \ponv \big |^2 \right) ^\h  \left( \sv  \big ( \pvnmx - \pvnmu \big )^2 \right) ^\h  \text{by Cauchy-Schwartz inequality} \\
  &= \frac{\norm{\pon}_2}{\e^\h} \left( \sv  \big ( \pvnmx - \pvnmu \big )^2 \e\right ) ^\h  \\
    &\leq  \frac{\norm{\pon}_2}{\e^\h} \left( \sv \sum_{w \in V} \big ( P_o^{n+m}(w)  - \ponmv \big )^2 P(w,v) \right ) ^\h \quad \text{by Lemma \ref{change_sum_over_all_vetices_} applied to $T=P^{n+m}$ } \\
      &\leq  \frac{\norm{\pon}_2}{\e^\h}  \bigg( 2  \inp{(I-P)\ponm}{\ponm} \bigg )^\h  \quad \text{by Lemma  \ref{change_u_v_sum_to_inner_product}}\ \\
      &= \left( \frac{2}{\e} \right)  ^\h   \norm{\pon}_2        \bigg(   \inp{(I-P)^\h\ponm}{(I-P)^\h\ponm} \bigg )^\h \quad \text{since $(I-P)^\h$ is self-adjoint} \\ 
       &= \left( \frac{2}{\e} \right)  ^\h      \norm{\pon}_2      \norm{ (I-P)^\h P^m \pon  }_2  \\ 
        &\leq  \left( \frac{2}{\e} \right)  ^\h   \norm{\pon}_2   \norm{\pon}_2      \norm{ (I-P)^\h P^m   }_{2\rightarrow 2}  \\ 
         &\leq  \norm{\pon}^2_2  \left( \frac{2}{\e} \right)  ^\h       \frac{1}{2 \sqrt{m}} \quad \text{ by Lemma \ref{bound_by_1_over_sqrt_2m} } \\
         &\leq  \frac{ \norm{\pon}^2_2  }{\sqrt{2m\e}}   =   \frac{ P^{2n}(o,o)  }{\sqrt{2m\e}} \quad \text{ by Lemma \ref{some_bounds} (\ref{2m_power_equal_2norm_squared_equal_inf_norm})}.
  \end{align*}

Let $o=x_0, x_1, \ldots , x_{|x|} = x$ be a shortest path from $o$ to $x$ with $x_i \sim x_{i+1}$. The above calculation tells us that
\begin{equation*}
|\pnnm(x) - \pnnm(o)| \leq \sum_{i=0}^{|x|-1} |\pnnm (x_{i+1}) - \pnnm (x_{i}) | \leq  |x|    \frac{ P^{2n}(o,o)  }{\sqrt{2m\e}}. 
\end{equation*}
Hence, for all $x \in V$, with $|x| \leq r(n,m)$, we have 
\[
|\pnnm(x) - \pnnm(o)|
\leq r(n,m)\,\frac{P^{2n}(o,o)}{\sqrt{2m\eta}}
= \left(\frac{\eta}{2} \right)^{\frac{1}{2}} m^{\frac{1}{2}} \frac{P^{2n+m}(o,o)}{P^{2n}(o,o)} \left(\frac{P^{2n}(o,o)}{\sqrt{2m\eta}}\right)
\leq \frac{P^{2n+m}(o,o)}{2}.
\]

It follows that, 
\[ \frac{P^{2n+m}(o,o)}{2} \leq \pnnm(x). \]  

Recall that there are $\beta(r(n,m))$ number of vertices  which are at most $r(n,m)$ away from $o$, so summing the inequality over the set of $\{x\in V \big| | x | \leq r(n,m) \}$ gives us
 
\begin{align*}
\sum_{\substack{{x \in V},\\  | x | \leq r(n,m) }} \frac{P^{2n+m}(o,o)}{2}   & \leq  \sum_{\substack{{x \in V},\\  | x | \leq r(n,m) }}  \pnnm(x)      \\
										\beta(r(n,m))	\frac{P^{2n+m}(o,o)}{2} 								& \leq  1 \\
									P^{2n+m}(o,o)	& \leq \frac{2}{\beta(r(n,m))}
\end{align*} 
which completes the proof.
\end{proof}

We will now find a bound for the decay of the sequence $P^{2m}(o,o)$ when a  vertex-transitive graph undergoes polynomial growth. The next theorem generalises \citep[Theorem~2.3]{MGA}, with the same constant $D$, and the proof follows theirs.

\begin{thm} \label{poly_growth_bound}
Let $(\G, P)$ be a symmetric transitive system. Recall that we defined $\eta=\inf\{P(u,v):u\sim v\}$. If the volume growth $\beta(j)$ satisfies 
\begin{equation} \label{volume_bound}
\beta(j) \geq aj^d, \quad \text{for all} \quad 1 \leq j \leq M 
\end{equation}

for some positive $a$, $d$ and $M$, then
\[P^{2m}(o,o) \leq \frac{D}{m^{d/2}} \quad \text{for all} \quad m \leq \frac{4M^2}{\e}\]
with $D = {2^{2+\frac{3}{2}d + \frac{d^2}{2}}}/(a \e ^{d/2}). $
\end{thm}

\begin{proof}

Let us first assume that $P$ is positive semi-definite. Set $A(n)=P^n(o,o)$. Then we have 
\begin{align}
A(2n+m) &\leq  \frac{2}{\beta(r(n,m))} \quad \text{by Lemma \ref{bound_on_p2nm} } \nonumber\\ 
&\leq \frac{2}{a} r(n,m)^{-d} \quad \text{using the volume bound (\ref{volume_bound}) }  \nonumber\\ 
&=\frac{2}{a}  \left \{ \left(\frac{\e}{2} \right)^{-1/2} m^{-1/2} \left( \frac{A(2n+m)}{A(2n)} \right) ^{-1}  \right\}^d, \label{bound_a2nm_initial}
\end{align}  
this holds for any $n$, $m$ with $r(n,m) \leq M$. From the definition for $r$ in Lemma \ref{bound_on_p2nm} and the bounds in Lemma \ref{some_bounds}, $r(n,m) \leq \left( \frac{\e m}{2} \right)^\h$, therefore, the formula (\ref{bound_a2nm_initial}) holds for $m \leq \frac{2M^2}{\e}$ and for all $n$.

Rewriting (\ref{bound_a2nm_initial}) to make  $A(2n+m)$  the subject, we get
\begin{equation} \label{repeat}
A(2n+m) \leq \left\{  \left(\frac{2}{a} \right)^\frac{1}{d} \left( \frac{2}{\e } \right)^\h \left( \frac{1}{ m} \right)^\h A(2n) \right\}^\theta,
\end{equation}
where $\theta = \frac{d}{1+d}$. Fix $m_0 \leq \dfrac{4M^2}{\eta}$. Let $k$ be such that $2^k \leq m_0 < 2^{k+1}$. Then, by using (\ref{repeat}) repeatedly, we have 

\begin{align*}
A(2^k) 
&= A(2^{k-1} + 2^{k-1} ) \\
&\leq \left\{  \left(\frac{2}{a} \right)^\frac{1}{d} \left( \frac{2}{\e } \right)^\h \left( \frac{1}{ 2^{k-1}} \right) ^\h A(2^{k-1}) \right\}^\theta \\
&\leq \left\{  \left(\frac{2}{a} \right)^\frac{1}{d}  \left( \frac{2}{\e }  \right)^\h \right\}^\theta \left\{  \frac{1}{ 2^{k-1}}   \right\}^\frac{\t}{2}  \left\{  \left(\frac{2}{a} \right)^\frac{1}{d} \left( \frac{2}{\e } \right)^\h \left( \frac{1}{ 2^{k-2}} \right)^\h A(2^{k-2}) \right\}^{\theta^2} \\
&\vdots\\
&\leq \left\{  \left(\frac{2}{a} \right)^\frac{1}{d}  \left( \frac{2}{\e }  \right)^\h \right\}^{S_1(\theta)} \left\{   2^{S_2(\theta)}   \right\}   A(2)^{\theta^{k-1}} \\
\end{align*}
where \[S_1(\theta)= \t + \t^2 + \ldots+ \t^{k-1} = \sum_{i=1}^{k-1} \t^i = d(1-\t^{k-1}) \leq d;\]  and \[S_2(\theta)= -\left[\frac{k-1}{2}\t + \frac{k-2}{2}\t^2+ \ldots \h \t^{k-1}\right] = \h \left( \sum_{i=1}^{k-1} i\t^i -k \sum_{i=1}^{k-1} \t^i \right) . \]

By differentiating $S_1(\theta)$ and multiplying the derivative by $\t$, we get that  \[\sum_{i=1}^{k-1} i\t^i = (d+d^2)\big( 1- k\theta^{k-1} + (k-1)\theta^k \big)  . \]

We will use this to show that  $S_2(\theta) \leq \frac{d +d^2 -kd  }{2}$ as follows:

\begin{align*}
2S_2(\theta) - \bigg( d + d^2 - kd \bigg)
&= (d+d^2)\big( 1- k\theta^{k-1} + (k-1)\theta^k \big) - k d(1-\theta^{k-1}) \\
&\quad - \bigg( d + d^2 - kd \bigg)\\
 &=(d +d^2) \big[  -k \t^{k-1} + (k-1)\t^k  \big] +kd \t^{k-1}  \\ 
 &=  -dk  \t^{k-1} + d(k-1)\t^k  - d^2 k \t^{k-1} + d^2 (k-1)\t^k  + kd \t^{k-1} \\
  &=   d(k-1)\t^k  - d^2 k \t^{k-1} + d^2 (k-1)\t^k   \\
 &= \bigg( d(k-1)\t  - d^2 k  + d^2 (k-1)\t  \bigg) \t^{k-1} \\
 &=  \bigg( (d+1)d(k-1)\t  - d^2 k   \bigg) \t^{k-1} \\
 &= \bigg( d^2(k-1)  - d^2 k   \bigg) \t^{k-1}   \qquad \text{since $\theta = \frac{d}{1+d}$} \\
 &= - d^2 \t^{k-1}  \leq 0.
\end{align*}

 We obtain the following bound when $P$ is positive semi-definite:
\begin{align*}
P^{m_0}(o,o)   &\leq  A(2^{k}) \qquad \text{by Lemma \ref{some_bounds} (\ref{2m+2_power_bdd_by_2m_power}),} \\
& \leq\left\{  \left(\frac{2}{a} \right)^\frac{1}{d}  \left( \frac{2}{\e }  \right)^\h \right\}^d \left\{   2^{\frac{d +d^2 -kd }{2}}   \right\} \\
&=  \left(\frac{2}{a} \right) \left( \frac{2}{\e }  \right)^\frac{d}{2}     2^{\frac{d +d^2}{2}}  \frac{1}{2^{\frac{kd }{2}} }  \frac{2^{\frac{d}{2}}}{2^{\frac{d}{2}}} \\
&=\frac{2^{1+\frac{d}{2}+\frac{d+d^2}{2}+\frac{d}{2}}}{a \e^\frac{d}{2}} \frac{1}{2^{\frac{(k+1)d }{2}} } \\
&\leq \frac{2^{1+ (3d+d^2)/2 }}{a (\e m_0)^{\frac{d}{2}}}  & \text{since $2^{k+1} \geq m_0$}.
\end{align*}

The above bound holds for the case when $P$ is a positive semi-definite. For general $P$,  Lemma \ref{P_plus} tells us that 
\begin{align*}
P^{2m}(o,o) &\leq   2  (P^+)^{2m}(o,o)  \\
& \leq 2 \frac{2^{1+ (3d+d^2)/2 }}{a \big ((\e/2) 2m \big)^{\frac{d}{2}}} \quad \text{ for } 2m \leq \frac{4M^2}{\e/2} \\
& \leq  \frac{2^{2+ (3d+d^2)/2 }}{a (\e m)^{\frac{d}{2}}} \quad \text{ for } m \leq \frac{4M^2}{\e}
\end{align*}

which completes the proof.

\end{proof}


\section{Bounds for the eigenvalues}
In this section, we will try to find two quantitative bounds for eigenvalues $\pi_2$ and $\pi_N$.
\subsection{Lower bound for the smallest eigenvalue}

\begin{lem} \label{pi_n-1}
Let $P$ be a symmetric stochastic matrix with rows and columns indexed by elements in $V$. Let  $1 = \pi_1   \geq \pi_2 \geq \ldots \geq \pi_{N}$ be the eigenvalues of $P$. Let $\delta = \inf\{P(u,u) | u \in V\}$, the minimum value of the diagonal elements of $P$. Then, the  smallest eigenvalue of $P$  satisfies
\[ \pi_{N} \geq -1 +2\delta.
\]
\end{lem}
\Cref{pi_n-1} is the second statement of \citep[Lemma~2.2]{MGA}; it holds for any symmetric stochastic matrix, so no group structure is needed.

\begin{proof}
If $\delta = 1$, then $P $ is the identity matrix with every eigenvalue being $1$, so the result clearly holds. If $\delta < 1$, let  $Q = (1-\delta)^{-1} (P - \delta I )$. Then $Q$ is also a symmetric and stochastic matrix. Hence the eigenvalues of $Q$ are bounded below by $-1$. 
The eigenvalues of $Q$ are $(1-\delta)^{-1} (\pi_i -\delta )$, for $1 \leq i \leq N$. In particular, we have
\[ (1-\delta)^{-1} (\pi_{N}-\delta) \geq -1
\]
which gives us the result. 
\end{proof}

\subsection{Lower bound and upper bound for the second biggest eigenvalue}

By looking at the ratio of Dirichlet forms, we obtain the following upper and lower bounds for the second biggest eigenvalue.

\begin{lem} \label{lem:bounds_on_Dirichlet_form_implies_bounds_on_eigenvalue}
Let $\G$ be a finite graph with order $N$. Let $P$ be the matrix associated with a symmetric random walk, with eigenvalue  $1 =\pi_1   \geq \pi_2 \geq \ldots \geq \pi_{N} $. Then, 
\begin{enumerate}[(i)]
\item \label{lower} For any non-constant real function $g$ on $V$, we have
\[1- \pi_2 \leq \frac{\mathscr{E} _P (g,g) }{\mathscr{E} _U (g,g)},\]
with equality when $g$ is an eigenfunction for $\pi_2$; in particular, $1-\pi_2$ is the minimum of this ratio over non-constant $g$.
\item \label{upper} Suppose $ \mathscr{E}  _U \leq A \mathscr{E}  _P $ for some constant $A$, then $\pi_i  \leq 1 - 1/A $, for all $i \geq 2$.
\end{enumerate}

\end{lem}

\begin{proof}
Part (i) is \cite[Lemma~13.7, Remark~13.8]{LPW}. For (ii), if $\mathscr{E}_U \leq A\,\mathscr{E}_P$, then $\mathscr{E}_P/\mathscr{E}_U \geq 1/A$ on non-constant functions, so (i) gives $1-\pi_2 \geq 1/A$, and $\pi_i \leq \pi_2$ for $i \geq 2$. See also \cite[(2.3)]{DSC4}.
\end{proof}

We will now find a bound for the ratio of Dirichlet form on a vertex-transitive graph. The following theorem generalises the diameter bound in \citep[Lemma~2.2]{MGA}, which is stated there for Cayley graphs. As in \cref{change_sum_over_all_vetices_}, the difference from the Cayley graph case is that we average over the stabiliser (the constant $C$ below).

\begin{thm} \label{thm:bounds_on_Dirichlet_form} 
Let $\G$ be a finite connected vertex-transitive graph with $\diam(\G)=\gamma$. Let $P$ be the transition matrix for a symmetric random walk with  $\eta=\inf\{P(u,v):u\sim v\}>0$ . Then we have 
\[\mathscr{E}  _U \leq \frac{ \gamma^2}{\e} \mathscr{E}  _P. \]
\end{thm}

\begin{proof}

 Since $\Aut(\G)$  acts transitively on $V$,  the size of the stabiliser of every $v$ in $V$ under the action $\Aut(\G)$ is the same, we denote $C$ to be this constant.

Take two vertices $x$, $y$ in $V$; suppose $x_0 = x$, $x_1$, $x_2$, $\ldots$ , $x_k = y$ is a shortest path from $x$ to $y$. 
Given $\alpha \in \Aut(\G)$, note that $\alpha(x_0)$, $\alpha(x_1)$, $\alpha(x_2)$, $\ldots$, $\alpha(x_k)$ is still a path, i.e. $\alpha(x_i) \sim \alpha(x_{i+1})$ for all $i$.

For any $\alpha \in \Aut(\G)$, we have 
\begin{align*} 
f(\alpha(x))-f(\alpha(y))  &= \big[f(\alpha(x_0))-f(\alpha(x_1)) \big] + \\
                            & \big[f(\alpha(x_1))-f(\alpha(x_2)) \big] +  \\
                            &\ldots\\
                             &  +   \big[f(\alpha(x_{k-1}))-f(\alpha(x_k)) \big]. 
\end{align*}
Then, squaring both sides and applying the Cauchy-Schwarz inequality, we obtain 
\begin{align} 
\big[ f(\alpha(x))-f(\alpha(y)) \big] ^2  \leq |k| \bigg\{   &\big[f(\alpha(x_0))-f(\alpha(x_1)) \big]^2 + \nonumber \\
                            & \big[f(\alpha(x_1))-f(\alpha(x_2)) \big]^2 +  \nonumber \\
                            &\ldots \nonumber \\
                             &  +   \big[f(\alpha(x_{k-1}))-f(\alpha(x_k)) \big] ^ 2 \bigg \}\label{path_inequality}. 
\end{align}
We will do some manipulations on both sides of the inequality (\ref{path_inequality}).

(\textbf{Step 1})We will firstly sum over all $\alpha \in \Aut(\G)$. On the left-hand side, we have

\begin{equation*} 
\sum_{\alpha \in \Aut(\G)} \big[ f(\alpha(x))-f(\alpha(y)) \big] ^2 = \sum_{u \in V}    \sum_{\substack{{\alpha \in \Aut(\G)},\\
                              \alpha(x)=u}}
\big[ f(u)-f(\alpha(y)) \big] ^2  
\end{equation*}

On the right-hand side, we have 
\begin{align*} 
 \sum_{\alpha \in \Aut(\G)}  |k| \bigg\{ \sum_{i=0}^{k-1}    \big[f(\alpha(x_i))-f(\alpha(x_{i+1})) \big]^2 \bigg\} 
                          &= |k|  \sum_{i=0}^{k-1} \sum_{u \in V}\sum_{\substack{{\alpha \in \Aut(\G)},\\
                              \alpha(x_i)=u}}   \big[ f(u)-f(\alpha(x_{i+1})) \big] ^2     \\
                              &\leq  |k|  \sum_{i=0}^{k-1} \sum_{u \in V} \sum_{\substack{{\alpha \in \Aut(\G)},\\
                              \alpha(x_i)=u}} \sum_{\substack{{v \in V},\\
                              v \sim u} }  \big[ f(u)-f(v) \big] ^2     \\
                              &= |k|^2 C \sum_{u \in V}  \sum_{\substack{{v \in V},\\
                              v \sim u} }  \big[ f(u)-f(v) \big] ^2\\                   
                              &\leq \gamma^2 C \sum_{u \in V}  \sum_{\substack{{v \in V},\\
                              v \sim u} }  \big[ f(u)-f(v) \big] ^2. 
\end{align*}

(\textbf{Step 2})Then we multiply both sides by $U(x,y) = \frac{1}{|V|}$ and sum over $y \in \G$. On the left-hand side, we have 
\begin{align*} 
 &\quad \sum_{y \in V} \sum_{u \in V}    \sum_{\substack{{\alpha \in \Aut(\G)},\\
                              \alpha(x)=u}} 
\big[ f(u)-f(\alpha(y)) \big] ^2  \frac{1}{|V|} \\
&=  \sum_{u \in V}    \sum_{\substack{{\alpha \in \Aut(\G)},\\  \alpha(x)=u}} \sum_{y \in V} \big[ f(u)-f(\alpha(y)) \big] ^2  \frac{1}{|V|} \\
&=  \sum_{u \in V}    \sum_{\substack{{\alpha \in \Aut(\G)},\\  \alpha(x)=u}} \sum_{v \in V} \big[ f(u)-f(v) \big] ^2  \frac{1}{|V|} &&\text{since $\alpha$ is an automorphism} \\
&= C \sum_{u \in V}   \sum_{v \in V} \big[ f(u)-f(v) \big] ^2  \frac{1}{|V|}. \\
\end{align*}

Note that there are no $y$ terms on the right-hand side,  multiplying by $ \frac{1}{|V|}$ and summing over $y \in \G$ has no effect:
\begin{align*} 
 \gamma^2 C \sum_{y \in V} \sum_{u \in V}  \sum_{\substack{{v \in V},\\ v \sim u} }  \big[ f(u)-f(v) \big] ^2   \frac{1}{|V|}
 &= \gamma^2 C \sum_{u \in V} \sum_{\substack{{v \in V},\\ v \sim u} }  \big[ f(u)-f(v) \big] ^2           \\        &\leq  \gamma^2 C \frac{1}{\e} \sum_{u \in V}   \sum_{\substack{{v \in V},\\ v \sim u} }  \big[ f(u)-f(v) \big] ^2  P(u,v). 
\end{align*}

(\textbf{Step 3})Finally, we divide both sides by $2C$; the left-hand side becomes  $ \mathscr{E}_U $, while the right-hand side gives us
\[
\gamma^2 \frac{1}{\eta} \bigg\{ \frac{1}{2} \sum_{u \in V} \sum_{\substack{v \in V\\ v \sim u}} [f(u)-f(v)]^2 P(u,v) \bigg\}
\leq \frac{1}{\eta} \gamma^2 \bigg\{ \frac{1}{2} \sum_{u \in V} \sum_{v \in V} [f(u)-f(v)]^2 P(u,v) \bigg\}
\leq \frac{\gamma^2}{\eta} \mathscr{E}_P.
\]

\end{proof}

\begin{cor} \label{pi1}
Let $\G$ be a finite connected vertex-transitive graph. Let $P$ be the matrix associated with a symmetric random walk. Suppose $\eta=\inf\{P(u,v):u\sim v\}>0$. Then the second biggest eigenvalue satisfies
\[ \pi_2 \leq  1-\frac{\e}{ \gamma^2}.
\]
\end{cor}

\begin{proof}
This follows from Lemma \ref{lem:bounds_on_Dirichlet_form_implies_bounds_on_eigenvalue} (\ref{upper}) and Theorem \ref{thm:bounds_on_Dirichlet_form}. 
\end{proof}

\section{Mixing time under Moderate growth}

\subsection{Proof of the main theorem}
We begin with some observations on the moderate growth property defined in \cref{sec:intro}. Clearly, any vertex-transitive graph $\G$ has $(A,d)$ moderate growth for some constants $A$ and $d$, in fact, one can just take $A = |\G|$. However, it turns out that many natural families of  vertex-transitive graphs have moderate growth with the same constants $A$ and $d$.

\begin{eg}
\begin{enumerate}
    \item For $m \in \mathbb{N}$, consider the Cayley graph $\Cayley(\mathbb{Z}_m, \{-1,0,1\})$. We have $\beta(j) = 2j+1$ and $\diam(\G) = m/2$. It follows that $\Cayley(\mathbb{Z}_m, \{-1,0,1\})$ has $(1,1)$-moderate growth.
    
    \item More generally, the Cayley graphs of finite nilpotent groups have $(A,d)$-moderate growth, where $A$ and $d$ depend only on the number of generators and the nilpotency class. This follows from \citep[Lemma~5.1 and Theorem~5.2]{MGA}.
    
    \item In the same paper \citep[Section~7.A]{MGA}, Diaconis and Saloff-Coste showed that the Cayley graphs of finite affine groups over $\mathbb{F}_p$ exhibit $(1,2)$-moderate growth, even though these groups do not have polynomial growth uniformly in $p$: indeed, the volume satisfies $\beta(j) \geq e^{cj}$ for $1 \leq j \leq \log p$, for a universal constant $c>0$ \citep[Section~7.B, Remark~3]{MGA}.
\end{enumerate}
\end{eg}

\cref{main_thm} tells us that, if the graph undergoes moderate growth, by putting some restrictions when travelling between adjacent   vertices, we can study the rate of convergence of the symmetric random walk to a uniform distribution. The proof below follows that of \citep[Theorem~3.1]{MGA}, with the ingredients above (\cref{poly_growth_bound} and \cref{pi1}) in place of their Cayley graph counterparts (Theorem~2.3 and Lemma~2.2 there). The lower bound uses the same eigenfunction argument and the same test function $g(x)=|x|$ as in their proof.

We present the proof here:
\begin{proof} [Proof of Theorem \ref{main_thm}]
\mbox{}

\begin{enumerate}[(i)]
\item
Combining Lemma \ref{pi_n-1} and Corollary \ref{pi1}, we obtain that $\pi_* \leq 1-  \frac{\e}{\gamma^2}$. 
Since $\G$ has moderate growth, $\beta(j) \geq \frac{N}{A\g^d} j^d$ for $1 \leq j \leq \g$.Substitute $M= \g$, $a = \frac{N}{A \g ^d }$     in Theorem \ref{poly_growth_bound} and take $ m= \frac{4\g^2}{\e}$, we obtain that 
\begin{equation*}
P^{2m}(o,o) \leq  \frac{2^{2+\frac{3}{2}d + \frac{d^2}{2}}}{          \left(\frac{N}{A \g ^d }\right)      \left(\e\frac{4\g^2}{\e} \right)^{d/2}} = 2 ^{2+ \frac{d+d^2}{2}} \frac{A}{N}. 
\end{equation*}
Recall an analysis result that the function $h:\R_{\geq 1} \rightarrow \R$, $h(x) = (1-\frac{1}{x})^x$ is an increasing function with $\lim_{x \rightarrow \infty} h(x) = e^{-1}$. Let $n$ be as in the statement of the theorem, and recall $c := \frac{n\eta}{\gamma^2}-4 > 0$. Set
\[
n' := n - m = \frac{c\gamma^2}{\eta},
\]
so that $\frac{n'\eta}{\gamma^2} = c$ and $n'+m=n$. Applying the bound (\ref{bound_aim}) with parameters $n'$ (for the eigenvalue term) and $m$ (for the return-probability term), we obtain
\begin{align*}
\|(P^{n}-U)_o\|^2_{2 } = \|(P^{n'+m}-U)_o\|^2_{2 } &\leq    2 ^{2+ \frac{d+d^2}{2}} \frac{A}{N} \left( 1-  \frac{\e}{\gamma^2}\right)^{2n'} \\
& = 2 ^{2+ \frac{d+d^2}{2}} \frac{A}{N} \left[ \left( 1-  \frac{\e}{\gamma^2}\right)^{\g^2 / \e}  \right] ^{2n'\e/\g^2} \\
& = 2 ^{2+ \frac{d+d^2}{2}} \frac{A}{N} \left[ \left( 1-  \frac{\e}{\gamma^2}\right)^{\g^2 / \e}  \right] ^{2c} \\
& \leq  2 ^{2+ \frac{d+d^2}{2}} \frac{A}{N} e ^{-2c}.
\end{align*}
Finally, taking the square root on both sides and noting that $\| U_o\|_{2} = N^{-\h}$ gives the desired result.

\item
Using Lemma \ref{holder_conjugate_norm_equality} (\ref{s_infty}), we have 
\begin{equation} \label{bound_1_inf}
\| (P^{k}-U)_o\|_{1} = \|P^k-U\|_{\infty \rightarrow \infty} = \max_{\normi{f} \leq 1} \| P^{k}(f)-U(f)\| _{ \infty}.
\end{equation}
Let $f$ be a real normalised eigenfunction for $P$ corresponding to   $\pi_2$ with respect to infinity norm, i.e. $Pf =  \pi_2f$ and $\normi{f} =1$. From Lemma \ref{eignvalue_of_p-U} (\ref{eignvalue_of_U}), we know that $U(f) =0$. Hence,  (\ref{bound_1_inf}) tells us that
\begin{equation} \label{tv_bound}
  \| (P^{k}-U)_o\|_{1} \geq  \| P^{k}(f)\| _{ \infty}  =  \pi_2^k . 
\end{equation}
Next, we will find a lower bound for $\pi_2$. We will show that 
\begin{equation} \label{pi1_bound}
\pi_2 \geq 1- \frac{D}{\g^2} \quad \text{with} \quad D= 4 ^{2d+1} A^2. 
\end{equation}

Let $g$ be the graph distance function on $V$, i.e. $g(v) = |v|$ for $v \in V$. Lemma \ref{lem:bounds_on_Dirichlet_form_implies_bounds_on_eigenvalue} (\ref{lower}) tells us that 
\[1- \pi_2 \leq \frac{\mathscr{E} _P (g,g) }{\mathscr{E} _U (g,g)}. \]
Recall that $\mathscr{E} _P (g,g) = \frac{1}{2} \sum_{u,v \in V} (g(u)-g(v))^2 P(u,v) $. 
We assumed that $P(u,v)=0$ for  any two non-adjacent vertices $u$, $v$ in $V$, it follows that $(g(u)-g(v))^2 \leq 1$ for $P(u,v)>0$. 
Hence,  $\mathscr{E} _P (g,g)  \leq  \frac{1}{2} \sum_{u,v \in V} 1 \cdot P(u,v) = \frac{N}{2}. $
On the other hand, 
\[\mathscr{E} _U (g,g)  = \frac{1}{2} \sum_{u,v \in V} (g(u)-g(v))^2 U(u,v) = \frac{1}{2N} \sum_{0 \leq i, j \leq \g} (i-j)^2 \partial B(i) \partial B(j), \]
where we write $\partial B(i)= |\{ v \in V : |v| = i\}|$, the number of vertices that are distance $i$ away from $o$. Consider $B_1 = \{ u \in V : |u| \leq \frac{\g}{4} \}$, and $B_2 = \{ v \in V :\frac{3\g}{4} \leq|v| \leq \g   \}$, then for any $u \in B_1$, and $v \in B_2$, we have $\big| |u|-|v| \big| \geq \frac{\g}{2}$. Restricting the double sum over $(i,j)$ to pairs coming from $u \in B_1$, $v \in B_2$ (and using the factor $\frac{1}{N}$ already present in $\mathscr{E}_U$), we obtain
\[
\mathscr{E} _U (g,g) \;\geq\; \frac{1}{N} \cdot \h \cdot \frac{\g^2}{4} |B_1||B_2| \;=\; \frac{\g^2}{8N} |B_1||B_2| .
\]
The moderate growth implies that 
\[ |B_1|=  \beta(\g/4)\geq \frac{(\g/4)^d N}{\g^d A} =A^{-1} 4^{-d} N. \]
Let $w \in V$ with $|w|= \g$. By triangle inequality, 
\[z \in B(w,\tfrac{\g}{4}) \Rightarrow d(o, z) \geq d(o, w) - d(w,z) \geq  \frac{3\g}{4} \Rightarrow  z \in  B_2, \]
so $B(w, \g/4) \subseteq B_2$. Hence, $|B_2|\geq |B(w, \g/4)|=\beta(\g/4) \geq A^{-1} 4^{-d} N$. Substituting these bounds,
\[
\mathscr{E} _U (g,g) \;\geq\; \frac{\g^2}{8N} \left( A^{-1} 4^{-d} N \right)^2 \;=\; \frac{\g^2 N}{2 \cdot 4^{2d+1} A^2},
\]
and therefore
\[
1 - \pi_2 \;\leq\; \frac{N/2}{\g^2 N / (2 \cdot 4^{2d+1} A^2)} \;=\; \frac{4^{2d+1} A^2}{\g^2},
\]
which proves (\ref{pi1_bound}). Then the inequality (\ref{tv_bound}) becomes 
\begin{align*}
\| (P^{n}-U)_o\|_{1 } \geq  \left( 1 - \frac{A^2 4^{2d+1}}{\g^2} \right)^n. 
\end{align*}
The assumption $\g \geq A 2^{2d+2}$ implies that $\frac{A^2 4^{2d+1}}{\g^2}  \leq \frac{1}{4}$. Furthermore, by looking at the Taylor expansion, we see that $1 - x \geq e^{-2x}$ for $0 \leq x \leq 1/2$, which concludes the proof.

\end{enumerate}

\end{proof}

\subsection{Mixing time}
For $\epsilon >0$ and $1 \leq p \leq \infty$, the \emph{$l^p$ mixing time} of system $(\G,P)$ is
\[\tau_p(\epsilon) = \min \{n \geq 0 :\| (P^{n}-U)_o\|_{p} \leq \epsilon \| U_o\|_{p}  \}.\]

Following from Theorem~\ref{main_thm}, we can deduce an upper bound and a lower bound for the mixing time for the random walk in the vertex-transitive graph.

\begin{cor} \label{eg_mixing_time}
\begin{enumerate} [(i)]
\item The mixing time for the system $(\G,P)$ described in Theorem \ref{main_thm} (\ref{main_upper}) is bounded above by \[\tau_2(\epsilon) \leq \frac{\gamma^2}{\eta} \log \frac{A^\h 2^{1+\frac{d+d^2}{4}}e^5} {\epsilon} .  \] 

\item The mixing time for the system $(\G,P)$ described in Theorem \ref{main_thm} (\ref{main_lower}) is bounded below by \[\tau_1(\epsilon) \geq \frac{\gamma^2}{2^{4d+3}A^2} \log \frac{1} {\epsilon}.  \] 
\end{enumerate}
\end{cor}

\begin{proof}
By Theorem \ref{main_thm}\eqref{main_upper}, $\|(P^n-U)_o\|_2/\|U_o\|_2 \le Be^{-c}$ with $c = \frac{n\eta}{\gamma^2}-4$, valid whenever $c>0$. Hence $\|(P^n-U)_o\|_2/\|U_o\|_2 \le \epsilon$ whenever
\[
Be^{-\left(\frac{n\eta}{\gamma^2}-4\right)} \le \epsilon
\quad\Longleftrightarrow\quad
n \ge n_0 := \frac{\gamma^2}{\eta}\left(\log\frac{B}{\epsilon}+4\right) = \frac{\gamma^2}{\eta}\log\frac{Be^4}{\epsilon}.
\]
Since $\tau_2(\epsilon)$ must be an integer, $\tau_2(\epsilon) = \lceil n_0 \rceil \le n_0+1$. As $\eta \le 1 \le \gamma^2$, we have $\gamma^2/\eta \ge 1$, so
\[
\tau_2(\epsilon) \le n_0 + \frac{\gamma^2}{\eta} = \frac{\gamma^2}{\eta}\log\frac{Be^4}{\epsilon} + \frac{\gamma^2}{\eta} = \frac{\gamma^2}{\eta}\log\frac{Be^5}{\epsilon},
\]
as claimed. Part (ii) follows the same way from Theorem \ref{main_thm}\eqref{main_lower}.
\end{proof}

We will use the following lemma to find an upper bound for the $l^\infty$ mixing time.

\begin{lem}
Let $(\G,P)$ be a symmetric transitive system and fix $o\in V$.\ Then:
\begin{enumerate} [(i)]
\item $\| (P^{n}-U)_o\|_{p}$ is non-increasing in $n$, for each $1\le p\le\infty$;
\item for each $\epsilon>0$, $\tau_p(\epsilon)$ is non-decreasing in $p$;
\item for all $n\ge 0$,
\[ \frac{\| (P^{2n}-U)_o\|_{\infty}}{\| U_o\|_{\infty} } = \left( \frac{\| (P^{n}-U)_o\|_{2}}{\| U_o\|_2 }  \right) ^2 ;\]
\item consequently, $\tau_\infty(\epsilon)\le 2\tau_2(\sqrt\epsilon)$ for every $\epsilon>0$.
\end{enumerate}
\end{lem}

\begin{proof}
\begin{enumerate} [(i)]
\item

According to inequality  \eqref{ssnorm_leq_1_norm}, we have
\begin{equation*}
\| (P^{n+1}-U)_o\|_{p}
= \| (P^{n}-U)_o P \|_{p}
\leq  \| (P^{n}-U)_o \|_{p}  \|P  \|_{p \rightarrow p}
\leq  \| (P^{n}-U)_o \|_{p}.
\end{equation*}
   
\item

The normalised norm $N^{-1/p}\|f\|_p=\big(\tfrac1N\sum_{u}|f(u)|^p\big)^{1/p}$ is non-decreasing in $p\in[1,\infty]$ by the power mean inequality (Jensen). Since $\|U_o\|_p=N^{1/p-1}$, we have $\|f\|_p/\|U_o\|_p=N\cdot N^{-1/p}\|f\|_p$. Applying this to $f=(P^n-U)_o$, the ratio $\|(P^n-U)_o\|_p/\|U_o\|_p$ is non-decreasing in $p$ for each $n$, so $\tau_p(\epsilon)$ is non-decreasing in $p$.
  
\item

For $n=0$ both sides equal $N-1$, so let $n\ge1$. For the numerator, Lemma \ref{pu_up_u} together with \eqref{signed_2m_bound} (applied with $T=P-U$ and $m=n$, using that $T^n=(P-U)^n=P^n-U$ satisfies \eqref{permutation_property}) tell us that 
\[\| (P^{2n}-U)_o\|_{\infty} = \|(P-U)^{2n}_o\|_{\infty} = \| (P-U)^{n}_o\|_2^2 = \| (P^{n}-U)_o\|_2^2.  \]
\item

Let $k=\tau_2(\sqrt{\epsilon})$. Then $\|(P^k-U)_o\|_2/\|U_o\|_2\le\sqrt{\epsilon}$, so by (iii) $\|(P^{2k}-U)_o\|_\infty/\|U_o\|_\infty\le\epsilon$, hence $\tau_\infty(\epsilon)\le 2k=2\tau_2(\sqrt{\epsilon})$.

\end{enumerate}
\end{proof}
Following from the above lemma, we can generalise our bounds in Corollary   \ref{eg_mixing_time} to the $l^p$ mixing time for $ 1 \leq p \leq \infty$,

\begin{cor} 
\begin{enumerate} [(i)]
\item The mixing time for the system $(\G,P)$ described in Theorem \ref{main_thm} (\ref{main_upper}) is bounded above by \[\tau_p(\epsilon) \leq \tau_\infty(\epsilon) \leq \frac{\gamma^2}{\eta} \log \frac{A\, 2^{2+\frac{d+d^2}{2}}\, e^{10}} {\epsilon} .   \] 

\item The mixing time for the system $(\G,P)$ described in Theorem \ref{main_thm} (\ref{main_lower}) is bounded below by \[\tau_p(\epsilon) \geq \tau_1(\epsilon)  \geq \frac{\gamma^2}{2^{4d+3}A^2} \log \frac{1} {\epsilon}.  \] 

\end{enumerate}
\end{cor}

\begin{proof}
For part (i), using $\tau_\infty(\epsilon) \leq 2\tau_2(\sqrt{\epsilon})$ together with Corollary \ref{eg_mixing_time}(i), we have
\[
\tau_\infty(\epsilon) \leq 2\cdot\frac{\gamma^2}{\eta}\log\frac{Be^5}{\sqrt{\epsilon}}
= \frac{\gamma^2}{\eta}\left(2\log\frac{Be^5}{\sqrt{\epsilon}}\right)
= \frac{\gamma^2}{\eta}\log\frac{(Be^5)^2}{\epsilon}
= \frac{\gamma^2}{\eta}\log\frac{B^2 e^{10}}{\epsilon}.
\]
Since $B^2 = A\,2^{2+\frac{d+d^2}{2}}$, this gives
\[
\tau_\infty(\epsilon) \leq \frac{\gamma^2}{\eta}\log\frac{A\,2^{2+\frac{d+d^2}{2}}\,e^{10}}{\epsilon},
\]
and combined with $\tau_p(\epsilon) \leq \tau_\infty(\epsilon)$, this proves part (i). Part (ii) follows the same way from Corollary \ref{eg_mixing_time}(ii) together with $\tau_p(\epsilon) \geq \tau_1(\epsilon)$.
\end{proof}

\section{Application: Vertex-transitive graph with large diameter}

We will use the following asymptotic notation: for $X$,$Y \in \R$, we say $X= \O (Y)$ if $X$ is at most a constant multiple of $Y$. If this constant depends on some other variable $\delta$, we indicate this with a subscript and write $X= \O_\delta (Y)$. We state the following result from \citep[Corollary~2.8]{FST}, which tells us that a vertex-transitive graph with large diameter has moderate growth.

\begin{thm} [large diameter implies moderate growth] \label{large_diam_thm}
Let $\G$ be a finite connected vertex-transitive graph. For every $\delta \geq 0$, there exists $n_0 = n_0(\delta)$ such that if $\diam(\G) \geq n_0$ and 

\begin{equation} \label{large_diam}
\diam(\G) \geq \left( \frac{|\G|}{\beta (1)} \right)^\delta
\end{equation}

then $\G$ has $( \O_\delta(1), \O_\delta(1) )$ - moderate growth.

\end{thm}

\begin{cor} [large diameter implies quadratic mixing time] \label{large_implies_quadratic}
Let $\delta \geq 0$. Let $(\G, P)$ be a symmetric transitive system, such that $\G$ also satisfies condition (\ref{large_diam}), then
\begin{enumerate} [(i)]
\item If $(\G,P)$ also satisfies the condition in Theorem \ref{main_thm} (\ref{main_upper}), then  $(\G,P)$ has quadratic mixing time, i.e. for $ 1 \leq p \leq \infty$,
\[\tau_p(\epsilon) =\O_\delta\left(  \frac{\g^2}{\eta} \log\O_\delta \left(\frac{1}{\epsilon} \right) \right ) \]
\item For a lower bound,  if $(\G,P)$ satisfies the condition in Theorem \ref{main_thm} (\ref{main_lower}), we have 
\[\tau_p(\epsilon) \geq C(\delta, \epsilon) \gamma^2 \]
where $C(\delta, \epsilon) = \frac{1}{2^{\O_\delta(1)} {\O_\delta(1)} } \log \frac{1}{\epsilon} = \frac{1}{{\O_\delta(1)} } \log \frac{1}{\epsilon} $.
\end{enumerate}

\end{cor}

\begin{rem}
At the moment, we are unable to make the bound in Corollary \ref{large_implies_quadratic} explicit. This is due to our reliance on the non-explicit Theorem~\ref{large_diam_thm}. In particular, we do not know how the bound for the mixing time depends on $\delta$. For a deeper discussion about this, we refer the reader to the remark in  \citep[p.~4]{FST}
\end{rem}

\footnotesize{
\bibliographystyle{abbrv}
\bibliography{Biblio}
}

\end{document}